\documentclass[preprint,11pt]{elsarticle}
\usepackage[margin=1in]{geometry}
\usepackage{amsmath,amssymb,amsthm}
\usepackage{enumitem}
\usepackage{tikz}
\usepackage{booktabs}
\usepackage{longtable}
\usepackage{hyperref}
\usetikzlibrary{arrows.meta,positioning,bending}
\hypersetup{colorlinks=true,linkcolor=blue!50!black,citecolor=blue!50!black,urlcolor=blue!50!black}

\numberwithin{equation}{section}
\newtheorem{theorem}[equation]{Theorem}
\newtheorem{proposition}[equation]{Proposition}
\newtheorem{lemma}[equation]{Lemma}
\newtheorem{corollary}[equation]{Corollary}
\newtheorem{definition}[equation]{Definition}
\newtheorem{remark}[equation]{Remark}
\newtheorem{example}[equation]{Example}

\newcommand{\RR}{\mathbb{R}}

\newcommand{\TT}{\mathcal{T}}
\newcommand{\WM}{Whitehead move}
\newcommand{\rank}{\operatorname{rank}}
\newcommand{\aff}{\operatorname{aff}}
\newcommand{\Aut}{\operatorname{Aut}}

\makeatletter
\def\ps@pprintTitle{%
  \let\@oddhead\@empty
  \let\@evenhead\@empty
  \def\@oddfoot{\reset@font\hfil\thepage\hfil}%
  \let\@evenfoot\@oddfoot}
\makeatother

\begin{document}

\begin{frontmatter}

\title{Circle-Pattern Linear Systems for the Twenty Trivalent
Delaunay Types of a Genus-Two Surface: Rank, Feasibility and
Polytope Dimension}

\author{Armando Rodado Amaris}
\ead{armando@rodadolab.com}
\affiliation{organization={Centre for Quantum Software and Information, University of Technology Sydney},
  city={Sydney},
  state={NSW},
  country={Australia}}

\begin{abstract}
The hyperelliptic quotient of a closed genus-two Riemann surface is a
hyperbolic sphere with six cone points of angle $\pi$. Its Delaunay
triangulation $T$ is dual to a trivalent Voronoi graph $G$ of one of
twenty combinatorial types. For each type we study the angle system
$Ax=\mathbf 1$, $x\ge0$, where $A$ is the vertex--edge incidence
matrix of $T$ and $\pi x_e$ is the angle subtended at a cone point by
the Voronoi edge dual to $e$, equivalently the intersection angle of
the two circumcircles adjacent across $e$. We prove that $\rank A=6$
for all twenty types, including the six whose triangulation has a
loop, so the affine solution space always has the dimension of
Teichm\"uller space. The polytope $P=\{Ax=\mathbf 1,\ x\ge0\}$ behaves
differently. It is six-dimensional for exactly ten types, the ten
Delaunay-realizable ones. For nine types it has dimension one, two or
three, because non-negativity forces some wedge angles to vanish, and
for the type of signature $(2,2,2,2,8,8)$ it is empty. We show that
five of the nine and the empty type never occur as Delaunay
triangulations, and that the other four occur only as triangulations
of non-trivalent Delaunay decompositions. Every vertex of every $P$
is half-integral, as for fractional perfect matching polytopes. For
the ten full-dimensional types every facet $x_e=0$ lies at an edge
where the Whitehead move is admissible, though not conversely, and the
move identifies the facets of the two polytopes; these identifications
pair the facet orbits of the ten cells under their automorphism
groups. The free unfolding of the gluing grows at rate about $7.093$.
All finite statements are checked by an independent program, and the
rank, infeasibility, forced coordinates, dimensions, vertices and
facets are verified in Lean~4.
\end{abstract}


\end{frontmatter}

\tableofcontents

\section{Introduction}

\subsection{Background}

The companion paper~\cite{Paper1} classifies the trivalent Delaunay
graphs arising from the six order-two cone points of the hyperelliptic
quotient $S=X/\iota$ of a closed genus-two Riemann surface $X$. Up to
embedded graph equivalence there are exactly two \emph{generic} types
(dual triangulation simple), connected by a Whitehead-move orbit, and
exactly twenty types once bigons are permitted at the trivalent
level~\cite[Theorem~8.1]{Paper1}. For the two generic types,
\cite[Section~6]{Paper1} attaches to the Delaunay triangulation $T$ a
linear system $A_ix=\mathbf 1$ in one variable per edge, the
normalized angle subtended at a cone point by the dual Voronoi edge,
and proves that it has rank six~\cite[Theorem~6.9]{Paper1}, the real
dimension of Teichm\"uller space $\TT_2$.

The present paper extends this construction to all twenty trivalent
types and asks what the rank computation does and does not tell us.
Rank is a statement about the affine subspace $\{x:A_ix=\mathbf1\}$
alone. The geometrically meaningful object is the \emph{polytope}
obtained by intersecting this subspace with the orthant $x\ge0$,
which is where the wedge angles of an actual Delaunay triangulation
lie. Our main finding is that the two questions have different
answers: the affine rank is six for all twenty types without
exception, while the polytope dimension varies sharply and, for one
type, the polytope is empty.

\subsection{Summary of results}

Section~\ref{sec:construction} sets up the angle system for a general
trivalent type, including the six types whose triangulation has a
loop (an edge of $G$ bordering the same face on both sides,
Remark~\ref{rem:loop}). Section~\ref{sec:rank} proves that
$\rank(A_i)=6$ for all twenty types (Theorem~\ref{thm:rankalways6}).
Section~\ref{sec:polytope} determines the dimension of each polytope
$P_i$: ten types have dimension six, nine have dimension one, two or
three, and one is infeasible (Theorem~\ref{thm:dimtable});
Section~\ref{sec:examples} works out one example of each regime, and
Section~\ref{sec:table} gives the full table and its geometric meaning
(Proposition~\ref{prop:degenerate}): the ten full-dimensional types
are exactly the ten Delaunay-realizable types of~\cite{Paper1}, six
types never occur as Delaunay triangulations at all, and four occur
only as triangulations refining a non-trivalent Delaunay
decomposition. Section~\ref{sec:vertices} shows that every vertex of
every $P_i$ is half-integral (Theorem~\ref{thm:halfintegral}) and
describes the vertices combinatorially. Section~\ref{sec:facets}
determines the facets of the ten full-dimensional polytopes, proves
that every facet lies at a Whitehead-admissible edge
(Theorem~\ref{thm:facetsubset}), and shows that a Whitehead move
identifies the corresponding facets of the two polytopes
(Theorem~\ref{thm:gluing}); vertices transfer along these
identifications, and all but eight vertices are generated from the
cube's (Theorem~\ref{thm:vertexgen}). Section~\ref{sec:unfolding}
shows that these identifications are compatible with the automorphism
groups, so the ten cells glue into a finite complex with orbifold
cells, and computes the growth rate of its free tree-like unfolding.
Section~\ref{sec:verification} describes the independent computational
check and the Lean~4 formalization, and Section~\ref{sec:discussion}
discusses open questions.

\subsection{Conventions}

The twenty types are numbered $0,\dots,19$ as in
Table~\ref{tab:edgelists}. This numbering extends the labels
$g_0=d_1$ (the cube, dual to the octahedron), $g_1=d_2$ (the stacked
type) and $g_8,g_{19}$ (the two types sharing the signature
$(2,2,3,3,7,7)$) of~\cite[Table~2]{Paper1}. The \emph{signature} of a
type is the sorted list of face degrees of $G$, equivalently the
vertex degrees of $T$, a loop counting twice.

\section{The angle system of a trivalent type}
\label{sec:construction}

Let $G$ be one of the twenty trivalent types and $T=G^*$ its dual: a
triangulation of the sphere with six vertices (the cone points, one
per face of $G$), twelve edges (one per edge of $G$) and eight
triangular faces (one per vertex of $G$). For the generic types $T$ is
simple; in general it may have multiple edges and loops.

\begin{definition}[Angle system]
\label{def:circlepattern}
Let $e$ be an edge of $T$ with endpoints $p,p'$ (possibly $p=p'$),
bordering the triangles $t,t'$ with circumcenters $O,O'$, which are
the endpoints of the Voronoi edge $e^*$ of $G$ dual to $e$. Set
\[
  x_e := \frac1\pi\cdot\angle_p(O,O'),
\]
the signed angle at $p$ subtended by $e^*$. Since the two
circumcircles of $t$ and $t'$ pass through $p$ with centers $O,O'$,
$\pi x_e$ is also their intersection angle at $p$, and Springborn's
interior angle is $\theta^*_e=\pi(1-x_e)$~\cite[Section~7.5]{Paper1}.
The \emph{angle system} of the type is
\[
  A_ix=\mathbf 1,\qquad x\ge0,
\]
where $A_i$ is the $6\times12$ vertex--edge incidence matrix of $T$, a
loop contributing $2$ to its vertex's row.
\end{definition}

The geometric content is that of~\cite[Proposition~6.7]{Paper1}. The
Voronoi cell of a cone point $p$ is star-shaped with respect to $p$,
and the geodesics from $p$ to the endpoints of its boundary edges cut
it into wedges, one for each end of an edge of $T$ at $p$; these
wedges fill the cone angle $\pi$ of~\cite[Proposition~2.2]{Paper1},
which is the row equation $\sum_{e\ni p}x_e=1$. For the Delaunay
triangulation every wedge angle is non-negative, and
by~\cite[Proposition~7.7]{Paper1} $x_e=0$ exactly when $O=O'$, that
is, when the four vertices of $t\cup t'$ are concyclic and the
Voronoi edge $e^*$ has length zero. So the wedge angles of a
Delaunay triangulation lie in the polytope
\[
  P_i := \{x\in\RR^{12} : A_ix=\mathbf1,\ x\ge0\},
\]
which is the polytope $D_i$ of~\cite[Section~6]{Paper1}; the bound
$x\le\mathbf 1$ used there is implied by the row equations. If the
Voronoi graph is trivalent, no four cone points of adjacent triangles
are concyclic, so $x$ lies in the relative interior $\{x>0\}$. A
Delaunay \emph{triangulation} of a non-trivalent Delaunay
decomposition (one obtained by adding diagonals to its non-triangular
cells) still gives a point of $P_i$, with $x_e=0$ exactly at the
added diagonals.

\begin{definition}[Reference labeling]
\label{def:canonlabel}
For each type we fix the labeling of Table~\ref{tab:edgelists}: the
vertices $1,\dots,8$ of $G$, its edges $e_1,\dots,e_{12}$ in
lexicographic order of their endpoints (parallel edges in a fixed
order), and the cone points $1,\dots,6$. All labeled statements below
(forced coordinates, vertex lists, worked examples) refer to this
labeling.
\end{definition}

The edge list of $G$ alone does not determine the type: three pairs of
types, $\{11,12\}$, $\{16,17\}$ and $\{18,19\}$, have isomorphic
underlying abstract graphs with different embeddings, as recorded
in~\cite[Section~8.1]{Paper1} and~\cite[\S7.1]{Rodado2007} (types
$18$ and $19$ even have literally the same edge list in our
labeling). For this reason Table~\ref{tab:edgelists} also lists, for
each edge $e_k$, the two cone points it joins, which is column $k$ of
$A_i$.

\begin{longtable}{@{}c c p{4.9in}@{}}
\toprule
type & signature & $e_k{=}uv/ab$: edge $uv$ of $G$, joining cone points $a,b$ of $T$ \\
\midrule
\endhead
0 & $(4,4,4,4,4,4)$ & {\scriptsize \texttt{e1{=}12/12, e2{=}13/13, e3{=}14/23, e4{=}25/14, e5{=}26/24, e6{=}35/15, e7{=}37/35, e8{=}46/26, e9{=}47/36, e10{=}58/45, e11{=}68/46, e12{=}78/56}} \\
1 & $(3,3,4,4,5,5)$ & {\scriptsize \texttt{e1{=}12/12, e2{=}13/23, e3{=}14/13, e4{=}25/14, e5{=}26/24, e6{=}34/35, e7{=}37/25, e8{=}48/15, e9{=}56/46, e10{=}58/16, e11{=}67/26, e12{=}78/56}} \\
2 & $(2,3,4,4,5,6)$ & {\scriptsize \texttt{e1{=}12/12, e2{=}13/13, e3{=}13/23, e4{=}24/24, e5{=}25/14, e6{=}36/12, e7{=}47/25, e8{=}48/45, e9{=}56/16, e10{=}58/46, e11{=}67/26, e12{=}78/56}} \\
3 & $(2,3,3,4,6,6)$ & {\scriptsize \texttt{e1{=}12/12, e2{=}13/13, e3{=}13/23, e4{=}24/14, e5{=}25/24, e6{=}36/12, e7{=}45/45, e8{=}47/15, e9{=}58/25, e10{=}67/16, e11{=}68/26, e12{=}78/56}} \\
4 & $(3,3,3,3,6,6)$ & {\scriptsize \texttt{e1{=}12/12, e2{=}13/23, e3{=}14/13, e4{=}25/24, e5{=}26/14, e6{=}34/35, e7{=}37/25, e8{=}47/15, e9{=}56/46, e10{=}58/26, e11{=}68/16, e12{=}78/12}} \\
5 & $(2,2,3,5,5,7)$ & {\scriptsize \texttt{e1{=}12/12, e2{=}13/13, e3{=}13/23, e4{=}24/24, e5{=}25/14, e6{=}36/12, e7{=}47/25, e8{=}47/45, e9{=}56/16, e10{=}58/46, e11{=}68/26, e12{=}78/24}} \\
6 & $(2,3,3,4,5,7)$ & {\scriptsize \texttt{e1{=}12/12, e2{=}13/13, e3{=}13/23, e4{=}24/24, e5{=}25/14, e6{=}35/12, e7{=}46/25, e8{=}47/45, e9{=}58/24, e10{=}67/56, e11{=}68/26, e12{=}78/46}} \\
7 & $(2,2,5,5,5,5)$ & {\scriptsize \texttt{e1{=}12/12, e2{=}13/13, e3{=}13/23, e4{=}24/14, e5{=}25/24, e6{=}36/12, e7{=}46/15, e8{=}47/45, e9{=}56/25, e10{=}58/45, e11{=}78/56, e12{=}78/46}} \\
8 & $(2,2,3,3,7,7)$ & {\scriptsize \texttt{e1{=}12/12, e2{=}13/13, e3{=}13/23, e4{=}24/14, e5{=}24/24, e6{=}35/12, e7{=}46/12, e8{=}57/25, e9{=}58/15, e10{=}67/26, e11{=}68/16, e12{=}78/56}} \\
9 & $(2,2,2,4,6,8)$ & {\scriptsize \texttt{e1{=}12/12, e2{=}13/23, e3{=}13/13, e4{=}24/14, e5{=}25/24, e6{=}36/12, e7{=}47/45, e8{=}47/15, e9{=}57/14, e10{=}58/12, e11{=}68/16, e12{=}68/26}} \\
10 & $(2,2,2,6,6,6)$ & {\scriptsize \texttt{e1{=}12/12, e2{=}13/13, e3{=}13/23, e4{=}24/14, e5{=}25/24, e6{=}36/12, e7{=}47/15, e8{=}47/45, e9{=}58/26, e10{=}58/46, e11{=}67/14, e12{=}68/24}} \\
11 & $(2,2,4,4,4,8)$ & {\scriptsize \texttt{e1{=}12/12, e2{=}13/13, e3{=}13/23, e4{=}24/24, e5{=}25/14, e6{=}35/12, e7{=}46/25, e8{=}47/45, e9{=}57/24, e10{=}68/26, e11{=}68/56, e12{=}78/25}} \\
12 & $(2,2,4,4,6,6)$ & {\scriptsize \texttt{e1{=}12/12, e2{=}13/13, e3{=}13/23, e4{=}24/24, e5{=}25/14, e6{=}35/12, e7{=}46/25, e8{=}47/45, e9{=}56/24, e10{=}68/45, e11{=}78/46, e12{=}78/56}} \\
13 & $(2,3,3,3,4,9)$ & {\scriptsize \texttt{e1{=}12/12, e2{=}13/13, e3{=}13/23, e4{=}23/12, e5{=}24/11, e6{=}45/14, e7{=}46/14, e8{=}57/15, e9{=}58/45, e10{=}67/16, e11{=}68/46, e12{=}78/56}} \\
14 & $(2,2,2,2,8,8)$ & {\scriptsize \texttt{e1{=}12/12, e2{=}13/13, e3{=}13/23, e4{=}24/14, e5{=}24/24, e6{=}35/12, e7{=}46/12, e8{=}57/25, e9{=}57/15, e10{=}68/16, e11{=}68/26, e12{=}78/12}} \\
15 & $(2,2,2,3,5,10)$ & {\scriptsize \texttt{e1{=}12/12, e2{=}13/13, e3{=}13/23, e4{=}24/12, e5{=}25/11, e6{=}36/12, e7{=}46/14, e8{=}46/24, e9{=}57/15, e10{=}58/15, e11{=}78/56, e12{=}78/16}} \\
16 & $(2,2,3,3,4,10)$ & {\scriptsize \texttt{e1{=}12/12, e2{=}13/13, e3{=}13/23, e4{=}24/14, e5{=}25/24, e6{=}35/12, e7{=}45/14, e8{=}46/11, e9{=}67/15, e10{=}68/15, e11{=}78/56, e12{=}78/16}} \\
17 & $(2,2,3,4,5,8)$ & {\scriptsize \texttt{e1{=}12/12, e2{=}13/23, e3{=}13/13, e4{=}24/24, e5{=}25/14, e6{=}34/12, e7{=}45/14, e8{=}56/44, e9{=}67/45, e10{=}68/45, e11{=}78/46, e12{=}78/56}} \\
18 & $(2,2,2,3,3,12)$ & {\scriptsize \texttt{e1{=}12/12, e2{=}13/23, e3{=}13/13, e4{=}23/12, e5{=}24/11, e6{=}45/14, e7{=}45/14, e8{=}56/11, e9{=}67/15, e10{=}68/15, e11{=}78/56, e12{=}78/16}} \\
19 & $(2,2,3,3,7,7)$ & {\scriptsize \texttt{e1{=}12/12, e2{=}13/13, e3{=}13/23, e4{=}23/12, e5{=}24/11, e6{=}45/14, e7{=}45/14, e8{=}56/44, e9{=}67/45, e10{=}68/45, e11{=}78/46, e12{=}78/56}} \\
\bottomrule
\caption{The reference labeling (Definition~\ref{def:canonlabel}).
For the edge $e_k$, the entry $uv/ab$ means that $e_k$ joins the
vertices $u,v$ of $G$ and separates the faces $a,b$ of $G$, so that
column $k$ of $A_i$ has its nonzero entries in rows $a$ and $b$ (a
single entry $2$ when $a=b$, a loop of $T$). Types $8$ and $19$ share
a signature, and types $18$ and $19$ share an edge list; the cone
points distinguish them.}
\label{tab:edgelists}
\end{longtable}

\begin{remark}[Multiple edges of $T$]
\label{rem:bigon}
If two faces of $G$ share two edges, $T$ has two edges with the same
endpoints and $A_i$ has two equal columns. Each edge has its own
Voronoi edge, kite and wedge angle, so the two variables are
independent. A cone point of degree two (a two-sided face of $G$,
bounded by two parallel edges of $G$) always comes with such a pair:
its two triangles share both of its edges and have two distinct third
edges with the same endpoints.
\end{remark}

\begin{remark}[Loops of $T$]
\label{rem:loop}
An edge $e^*$ of $G$ can border the same face on both sides (a bridge
of $G$), giving $T$ a loop $e$ at the corresponding cone point $p$.
Both ends of $e$ are at $p$, so the Voronoi cell of $p$ meets $e^*$
from both sides and contains two wedges subtended by it. They have the
same angle $\pi x_e$, because the kite of $e$ is symmetric under the
reflection in $e^*$, the perpendicular bisector of $e$. The
coefficient $2$ in $A_i$ records exactly these two wedges. Six of the
twenty types have loops (Table~\ref{tab:results}).
\end{remark}

\section{The rank is always six}
\label{sec:rank}

\begin{lemma}[Incidence-matrix rank, with loops]
\label{lem:rankwithloops}
Let $H$ be a connected graph, possibly with loops and multiple edges,
with unsigned vertex--edge incidence matrix $M$, a loop contributing
$2$ to its vertex's row. If $H$ has a loop or an odd cycle, then
$\rank(M)=|V(H)|$; if $H$ is loopless and bipartite,
$\rank(M)=|V(H)|-1$.
\end{lemma}

\begin{proof}
A left null vector $y$ satisfies $y_u+y_v=0$ for every edge $uv$ with
$u\ne v$, and $2y_v=0$ for every loop at $v$. If $H$ has a loop at
$v_0$, then $y_{v_0}=0$, and connectedness propagates $y\equiv0$
along the edges, so $M$ has full row rank. If $H$ is loopless, this
is~\cite[Lemma~6.8]{Paper1}, whose proof does not use simplicity.
\end{proof}

\begin{theorem}[$\rank(A_i)=6$ for all twenty types]
\label{thm:rankalways6}
For each of the twenty trivalent types, $\rank(A_i)=6$, so the affine
solution space $\{x:A_ix=\mathbf1\}$ has dimension $12-6=6$.
\end{theorem}

\begin{proof}
$T$ is connected, and the boundary of any of its faces is a closed
walk of length three. A closed walk of odd length contains a loop or
an odd cycle, whether or not $T$ is simple. Apply
Lemma~\ref{lem:rankwithloops} with $H=T$.
\end{proof}

This extends~\cite[Theorem~6.9]{Paper1} from the two generic types to
all twenty. We also verified it by exact rational computation of each
of the twenty ranks, and in Lean (Section~\ref{sec:verification}).

\section{Feasibility and the dimension of the polytope}
\label{sec:polytope}

Theorem~\ref{thm:rankalways6} determines the dimension of the affine
subspace $\{x:A_ix=\mathbf1\}$, but $\dim P_i=\dim\aff(P_i)$ can be
strictly smaller if non-negativity forces some coordinates to a
constant value on all of $P_i$, and $P_i$ can be empty. Both happen.

\begin{definition}[Forced coordinate]
A coordinate $x_e$ is \emph{forced} (to $c$) if $x_e=c$ for every
$x\in P_i$. The forced values that occur are $0$ and $1$.
\end{definition}

\begin{proposition}
\label{prop:dimformula}
If $P_i\ne\emptyset$ and $F$ is the set of forced coordinates, then
\[
  \dim P_i = 12 - \rank\big([\,A_i;\ (e_j)_{j\in F}\,]\big),
\]
where $e_j$ is the $j$-th standard basis row vector.
\end{proposition}

\begin{proof}
The affine hull of a polyhedron $\{x:A x=b,\ x\ge0\}$ is cut out by
the equations $Ax=b$ together with the inequalities that hold with
equality on the whole polyhedron (its implicit
equalities)~\cite[Section~8.2]{Schrijver1986}. Here these are the
equations $x_j=0$ for the coordinates forced to $0$; a coordinate
forced to $1$ at a cone point $p$ forces every other coordinate at
$p$ to $0$, so adding $x_j=1$ as well does not change the affine
space. The displayed formula is its dimension.
\end{proof}

We determined $F$ for each type by linear programming (maximizing and
minimizing each $x_e$ over $P_i$) and confirmed it by exact rational
arithmetic and in Lean, where each forced value is a theorem
(Section~\ref{sec:verification}).

\begin{theorem}[Polytope dimensions]
\label{thm:dimtable}
Of the twenty types, ten have $\dim P_i=6$; nine have
$\dim P_i\in\{1,2,3\}$, with the forced coordinates and values of
Table~\ref{tab:results}; and the type of signature $(2,2,2,2,8,8)$ has
$P_i=\emptyset$.
\end{theorem}

The proof is Sections~\ref{sec:examples}--\ref{sec:table}: three
worked derivations, one for each regime, and the full table obtained
by the same method.

\section{Worked examples}
\label{sec:examples}

We write $R_p$ for the row equation $\sum_{e\ni p}x_e=1$ of the cone
point $p$, in the labeling of Table~\ref{tab:edgelists}.

\subsection{A full-dimensional example: the cube}

For $G=d_1$, $T$ is the octahedron and every cone point has degree
four, so the point $x_e=\tfrac14$ for all twelve $e$ satisfies
$A_0x=\mathbf1$. It is strictly positive, so no coordinate is forced
and $\dim P_0=6$ by Theorem~\ref{thm:rankalways6}. It is the wedge
vector of the explicit octahedral configuration
of~\cite[Proposition~7.2]{Paper1}.

\subsection{A partially degenerate example}
\label{subsec:example3}

Consider type $3$, of signature $(2,3,3,4,6,6)$. Its cone point $3$
has degree two, with edges $e_2,e_3$ (the parallel edges $e_2,e_3$ of
$G$ bound the corresponding two-sided face). Adding the row equations
of the cone points $1,2,5$ and subtracting those of $3,4,6$ gives
\[
  (R_1+R_2+R_5)-(R_3+R_4+R_6):\qquad
  2\big(x_{e_1}+x_{e_6}+x_{e_8}+x_{e_9}\big)=0.
\]
Since all four terms are non-negative,
\[
  x_{e_1}=x_{e_6}=x_{e_8}=x_{e_9}=0
\]
on all of $P_3$. These are three independent equalities beyond
$A_3x=\mathbf1$ (the fourth follows from them and the row equations),
so $\dim P_3=12-(6+3)=3$ by Proposition~\ref{prop:dimformula}. A point
in the relative interior is
\[
  x_{e_1}=x_{e_6}=x_{e_8}=x_{e_9}=0,\qquad
  x_{e_2}=x_{e_3}=x_{e_7}=x_{e_{12}}=\tfrac12,\qquad
  x_{e_4}=x_{e_5}=x_{e_{10}}=x_{e_{11}}=\tfrac14,
\]
which satisfies $A_3x=\mathbf1$ and is positive off the forced
coordinates.

\subsection{The infeasible example}
\label{subsec:example14}

Consider type $14$, of signature $(2,2,2,2,8,8)$. Its cone points
$3,4,5,6$ have degree two, and the remaining eight edge-ends at the
cone points $1$ and $2$ include the four parallel edges
$e_1,e_6,e_7,e_{12}$ joining them. Adding the row equations of $1$
and $2$ and subtracting those of $3,4,5,6$ gives
\[
  (R_1+R_2)-(R_3+R_4+R_5+R_6):\qquad
  2\big(x_{e_1}+x_{e_6}+x_{e_7}+x_{e_{12}}\big)=2-4=-2 .
\]
The left side is non-negative on the orthant, so $P_{14}=\emptyset$.
In words: the four degree-two cone points use up four of the twelve
edges entirely, and the cone angles at $1$ and $2$ cannot then be
filled by non-negative wedges.

\begin{corollary}
\label{cor:norealiz}
No hyperbolic sphere with six cone points of angle $\pi$ has a
Delaunay triangulation of type $14$, even allowing concyclic
quadruples (a triangulation refining a non-trivalent Delaunay
decomposition).
\end{corollary}

\begin{proof}
The wedge angles of such a triangulation would be a point of
$P_{14}$ (Section~\ref{sec:construction}).
\end{proof}

\section{The full results table}
\label{sec:table}

\begin{longtable}{@{}c c c c c l@{}}
\toprule
type & signature & $\rank(A_i)$ & $\dim P_i$ & loops & forced coordinates \\
\midrule
\endhead
0  & $(4,4,4,4,4,4)$   & 6 & 6 & ---  & --- \\
1  & $(3,3,4,4,5,5)$   & 6 & 6 & ---  & --- \\
2  & $(2,3,4,4,5,6)$   & 6 & 6 & ---  & --- \\
3  & $(2,3,3,4,6,6)$   & 6 & 3 & ---  & $e_1,e_6,e_8,e_9=0$ \\
4  & $(3,3,3,3,6,6)$   & 6 & 6 & ---  & --- \\
5  & $(2,2,3,5,5,7)$   & 6 & 2 & ---  & $e_1,e_4,e_5,e_6,e_{12}=0$ \\
6  & $(2,3,3,4,5,7)$   & 6 & 6 & ---  & --- \\
7  & $(2,2,5,5,5,5)$   & 6 & 6 & ---  & --- \\
8  & $(2,2,3,3,7,7)$   & 6 & 1 & ---  & $e_1,e_6,e_7,e_8,e_9,e_{10},e_{11}=0$; $e_{12}=1$ \\
9  & $(2,2,2,4,6,8)$   & 6 & 1 & ---  & $e_1,e_4,e_5,e_6,e_8,e_9,e_{10}=0$; $e_7=1$ \\
10 & $(2,2,2,6,6,6)$   & 6 & 1 & ---  & $e_1,e_4,e_5,e_6,e_{11},e_{12}=0$ \\
11 & $(2,2,4,4,4,8)$   & 6 & 3 & ---  & $e_1,e_6,e_7,e_{12}=0$ \\
12 & $(2,2,4,4,6,6)$   & 6 & 6 & ---  & --- \\
13 & $(2,3,3,3,4,9)$   & 6 & 6 & $e_5$ & --- \\
14 & $(2,2,2,2,8,8)$   & 6 & --- & --- & \textbf{infeasible} \\
15 & $(2,2,2,3,5,10)$  & 6 & 1 & $e_5$ & $e_1,e_4,e_5,e_6,e_9,e_{10},e_{12}=0$; $e_{11}=1$ \\
16 & $(2,2,3,3,4,10)$  & 6 & 2 & $e_8$    & $e_1,e_6,e_8,e_9,e_{10},e_{12}=0$; $e_{11}=1$ \\
17 & $(2,2,3,4,5,8)$   & 6 & 6 & $e_8$    & --- \\
18 & $(2,2,2,3,3,12)$  & 6 & 1 & $e_5,e_8$ & $e_1,e_3,e_4,e_5,e_8,e_9,e_{10},e_{12}=0$; $e_2,e_{11}=1$ \\
19 & $(2,2,3,3,7,7)$   & 6 & 6 & $e_5,e_8$ & --- \\
\bottomrule
\caption{The angle system of each of the twenty types, in the labeling
of Table~\ref{tab:edgelists}: rank, polytope dimension, the loops of
$T$, and the coordinates that take the same value at every point of
$P_i$ (Proposition~\ref{prop:dimformula}).}
\label{tab:results}
\end{longtable}

Ten types have $\dim P_i=6$; the remaining nine feasible types have
$\dim P_i\in\{1,2,3\}$ (five of dimension $1$, two of dimension $2$,
two of dimension $3$); one type is infeasible. Every entry was
obtained twice, by exact rational elimination and by floating-point
linear programming, and is verified in Lean
(Section~\ref{sec:verification}).

The table has a direct geometric reading. Write $Z_i$ for the set of
coordinates forced to $0$.

\begin{proposition}[What the degenerate types are]
\label{prop:degenerate}
\leavevmode
\begin{enumerate}[label=(\alph*)]
\item The ten full-dimensional types $0,1,2,4,6,7,12,13,17,19$ are
exactly the types with a strictly positive solution of
$A_ix=\mathbf 1$, and exactly the ten Delaunay-realizable trivalent
types of~\cite[Section~7.5]{Paper1}. For each of the other ten types,
no hyperbolic sphere with six cone points of angle $\pi$ has a
trivalent Voronoi graph of that type.
\item For the types $8,9,15,16,18$, which have a coordinate forced to
$1$, and for the infeasible type $14$, no such sphere has a Delaunay
triangulation of that type at all, even one refining a non-trivalent
Delaunay decomposition.
\item For the types $3,5,10,11$, the graph $T_i\smallsetminus Z_i$
is connected, with all vertex degrees at least two, and $P_i$ is the
polytope of the contracted graph $G_i/Z_i$, which has a strictly
positive solution. Every Delaunay triangulation of one of these types
refines a Delaunay decomposition whose $1$-skeleton is contained in
$T_i\smallsetminus Z_i$; it has $12-|Z_i|=8,7,6,8$ edges.
\end{enumerate}
\end{proposition}

\begin{proof}
(a) $\dim P_i=6=\dim\{A_ix=\mathbf1\}$ holds exactly when no
coordinate is forced, and then the average of points of $P_i$
maximizing each coordinate in turn is strictly positive; conversely a
strictly positive solution forces nothing. A trivalent Voronoi graph
gives a strictly positive point of $P_i$ (Section~\ref{sec:construction}),
which rules out the ten lower-dimensional or empty types. The ten
full-dimensional types are the ten types found Delaunay realizable
in~\cite[Section~7.5]{Paper1} by Springborn's closing conditions,
and the same ten have strictly positive solutions in the independent
check of~\cite[Remark~8.12]{Paper1}.

(b) A Delaunay triangulation of type $i$ gives $x\in P_i$, and the
edges of its Delaunay decomposition are the edges with $x_e>0$. For
the five types with a coordinate forced to $1$, the graph
$T_i\smallsetminus Z_i$ is disconnected and has vertices of degree
one (checked directly from Table~\ref{tab:edgelists}). But the
$1$-skeleton of the Delaunay decomposition contains all six cone
points and is connected, and by~\cite[Theorem~6.5]{Paper1} the
Voronoi graph has no monogon, that is, $T$ minus the zero edges has no
vertex of degree one. Type $14$ is Corollary~\ref{cor:norealiz}.

(c) Deleting the edges $Z_i$ from $T_i$ is dual to contracting them
in $G_i$, and $P_i$ lies in $\{x_{Z_i}=0\}$, so dropping these
coordinates identifies $P_i$ with the polytope of
$T_i\smallsetminus Z_i$. Its relative interior is positive off
$Z_i$. The edges of the Delaunay decomposition are those with
$x_e>0$, a subset of the complement of $Z_i$.
\end{proof}

So the nine lower-dimensional polytopes are not new cells: four of
them are the closed polytopes of realizable types with fewer edges,
one of which, for type $10$, is the six-cycle dual to the graph
$\Theta_6$ at the bottom of the contraction
hierarchy~\cite[Theorem~8.6]{Paper1}, and the other five correspond
to no configuration at all.

\section{Vertices of the polytopes}
\label{sec:vertices}

The extreme points of the polytopes $P_i$ have a common structure,
which is that of the fractional perfect matching polytope of a graph,
here for the triangulation $T$ and extended to graphs with loops.

\begin{theorem}[Half-integral vertices]
\label{thm:halfintegral}
For every one of the twenty types, every vertex of $P_i$ has all
twelve coordinates in $\{0,\tfrac12,1\}$.
\end{theorem}

For simple graphs this is Balinski's theorem on the fractional
matching polytope~\cite{Balinski1965}; see
also~\cite[Theorem~30.2]{Schrijver2003}. We give the proof, which
needs only a small change for loops.

\begin{proof}
Let $x$ be a vertex of $P_i$, let $F=\{e:0<x_e<1\}$ be its fractional
support, and let $V_F$ be the set of cone points incident to $F$. At
each $p\in V_F$ the integral coordinates in the row equation sum to an
integer, so the fractional terms $\sum_{e\in F,\,e\ni p}x_e$ (a loop
counted twice) sum to an integer, which is $1$ since they lie
strictly between $0$ and $2$ and are positive. In particular no
$p\in V_F$ meets $F$ in a single non-loop edge, so every component of
the graph $(V_F,F)$ has at least as many edges as vertices.

Since $x$ is a vertex, the columns of $A_i$ indexed by $F$ are
linearly independent: otherwise $x$ could be moved in both directions
along a kernel vector supported on $F$. Restricted to the rows $V_F$
(the other rows vanish on these columns), this is the incidence
matrix of $(V_F,F)$, so by Lemma~\ref{lem:rankwithloops}, applied to
each component, every component has at most as many edges as
vertices, with equality only if it has a loop or an odd cycle. Hence
every component has exactly as many edges as vertices and minimum
degree two, so it is a single cycle (a single loop being a cycle of
length one), and it is not bipartite, so the cycle is odd.

On an odd cycle with fractional values $t_1,\dots,t_\ell$ in cyclic
order, the equations $t_j+t_{j+1}=1$ (indices mod $\ell$) have the
unique solution $t_j=\tfrac12$; on a loop, $2t_1=1$. Hence every
fractional coordinate equals $\tfrac12$.
\end{proof}

\begin{corollary}[Combinatorial description of the vertices]
\label{cor:vertexdescription}
The vertices of $P_i$ correspond bijectively to the spanning
subgraphs of $T$ each of whose components is a single edge, an odd
cycle of length at least three, or a single loop: the vertex is $1$ on
the single-edge components, $\tfrac12$ on the edges of the other
components, and $0$ elsewhere. If such a subgraph has $m$ single-edge
components and its other components cover $h$ cone points, the vertex
has $m$ coordinates equal to $1$ and $h$ equal to $\tfrac12$, and
$2m+h=6$.
\end{corollary}

\begin{proof}
The proof of Theorem~\ref{thm:halfintegral} shows that the support of
a vertex is such a subgraph. Conversely, the vector described is in
$P_i$, and its support has linearly independent columns by
Lemma~\ref{lem:rankwithloops} (each component is an edge, or has a
loop or an odd cycle and as many edges as vertices), so it is a basic
feasible solution, hence a vertex. Different subgraphs give different
supports, hence different vertices.
\end{proof}

Each odd component covers an odd number of cone points and $6$ is
even, so the number of odd components is even and
$h\in\{0,2,4,6\}$. For the fourteen loop-free types the odd components
are cycles of length at least three, so $h\in\{0,6\}$, with $h=6$ only
for two disjoint triangles. Table~\ref{tab:vertices} lists the
results.

\begin{table}[htbp]
\centering
\small
\begin{tabular}{@{}c c r l@{}}
\toprule
type & signature & $|V(P_i)|$ & $(m,h)$-distribution \\
\midrule
0  & $(4,4,4,4,4,4)$   & 12 & $(3,0){:}\,8$,\ \ $(0,6){:}\,4$ \\
1  & $(3,3,4,4,5,5)$   & 9  & $(3,0){:}\,7$,\ \ $(0,6){:}\,2$ \\
2  & $(2,3,4,4,5,6)$   & 7  & $(3,0){:}\,5$,\ \ $(0,6){:}\,2$ \\
3  & $(2,3,3,4,6,6)$   & 4  & $(3,0){:}\,4$ \\
4  & $(3,3,3,3,6,6)$   & 8  & $(3,0){:}\,6$,\ \ $(0,6){:}\,2$ \\
5  & $(2,2,3,5,5,7)$   & 3  & $(3,0){:}\,3$ \\
6  & $(2,3,3,4,5,7)$   & 7  & $(3,0){:}\,5$,\ \ $(0,6){:}\,2$ \\
7  & $(2,2,5,5,5,5)$   & 8  & $(3,0){:}\,4$,\ \ $(0,6){:}\,4$ \\
8  & $(2,2,3,3,7,7)$   & 2  & $(3,0){:}\,2$ \\
9  & $(2,2,2,4,6,8)$   & 2  & $(3,0){:}\,2$ \\
10 & $(2,2,2,6,6,6)$   & 2  & $(3,0){:}\,2$ \\
11 & $(2,2,4,4,4,8)$   & 4  & $(3,0){:}\,4$ \\
12 & $(2,2,4,4,6,6)$   & 8  & $(3,0){:}\,4$,\ \ $(0,6){:}\,4$ \\
13 & $(2,3,3,3,4,9)$   & 7  & $(3,0){:}\,4$,\ \ $(1,4){:}\,1$,\ \ $(0,6){:}\,2$ \\
14 & $(2,2,2,2,8,8)$   & 0  & --- \\
15 & $(2,2,2,3,5,10)$  & 2  & $(3,0){:}\,2$ \\
16 & $(2,2,3,3,4,10)$  & 3  & $(3,0){:}\,3$ \\
17 & $(2,2,3,4,5,8)$   & 9  & $(3,0){:}\,3$,\ \ $(1,4){:}\,2$,\ \ $(0,6){:}\,4$ \\
18 & $(2,2,2,3,3,12)$  & 2  & $(3,0){:}\,2$ \\
19 & $(2,2,3,3,7,7)$   & 11 & $(3,0){:}\,2$,\ \ $(1,4){:}\,4$,\ \ $(2,2){:}\,1$,\ \ $(0,6){:}\,4$ \\
\bottomrule
\end{tabular}
\caption{Vertex counts and $(m,h)$-distributions, computed by
enumerating all $\binom{12}{6}=924$ choices of six basic columns and
keeping the feasible basic solutions, and cross-checked for all twenty
types against the count of subgraphs in
Corollary~\ref{cor:vertexdescription}. The entries $(1,4)$, at types
$13,17,19$, come from a loop together with a disjoint triangle, and
the entry $(2,2)$, at type $19$, from two loops.}
\label{tab:vertices}
\end{table}

For the two generic types the vertices are as follows.

\begin{example}[The cube, type $0$]
\label{ex:cubevertices}
$\Aut^+(T)$, of order $24$, acts on the twelve vertices of $P_0$ with
two orbits: the eight perfect matchings of the octahedron, with
stabilizers of order $3$, and the four pairs of opposite faces, with
stabilizers of order $6$. In the labeling of
Table~\ref{tab:edgelists}, the twelve vertices are
\[
\arraycolsep=2.2pt
\begin{array}{cccccccccccc}
\smash{\tfrac12} & \smash{\tfrac12} & \smash{\tfrac12} & 0 & 0 & 0 & 0 & 0 & 0 & \smash{\tfrac12} & \smash{\tfrac12} & \smash{\tfrac12} \\
\smash{\tfrac12} & 0 & 0 & \smash{\tfrac12} & \smash{\tfrac12} & 0 & \smash{\tfrac12} & 0 & \smash{\tfrac12} & 0 & 0 & \smash{\tfrac12} \\
1 & 0 & 0 & 0 & 0 & 0 & 1 & 0 & 0 & 0 & 1 & 0 \\
1 & 0 & 0 & 0 & 0 & 0 & 0 & 0 & 1 & 1 & 0 & 0 \\
0 & \smash{\tfrac12} & 0 & 0 & \smash{\tfrac12} & \smash{\tfrac12} & \smash{\tfrac12} & \smash{\tfrac12} & 0 & 0 & \smash{\tfrac12} & 0 \\
0 & 1 & 0 & 0 & 1 & 0 & 0 & 0 & 0 & 0 & 0 & 1 \\
0 & 1 & 0 & 0 & 0 & 0 & 0 & 1 & 0 & 1 & 0 & 0 \\
0 & 0 & \smash{\tfrac12} & \smash{\tfrac12} & 0 & \smash{\tfrac12} & 0 & \smash{\tfrac12} & \smash{\tfrac12} & \smash{\tfrac12} & 0 & 0 \\
0 & 0 & 1 & 1 & 0 & 0 & 0 & 0 & 0 & 0 & 0 & 1 \\
0 & 0 & 1 & 0 & 0 & 1 & 0 & 0 & 0 & 0 & 1 & 0 \\
0 & 0 & 0 & 1 & 0 & 0 & 1 & 1 & 0 & 0 & 0 & 0 \\
0 & 0 & 0 & 0 & 1 & 1 & 0 & 0 & 1 & 0 & 0 & 0
\end{array}
\]
Their average is $(\tfrac14,\dots,\tfrac14)$, the symmetric point of
the octahedral configuration~\cite[Proposition~7.2]{Paper1}, as it
must be: the average is fixed by $\Aut^+(T)$, which acts transitively
on the edges.
\end{example}

\begin{example}[The stacked type, type $1$]
\label{ex:stackedvertices}
In the same way, the nine vertices of $P_1$ are
\[
\arraycolsep=2.2pt
\begin{array}{cccccccccccc}
1 & 0 & 0 & 0 & 0 & 1 & 0 & 0 & 1 & 0 & 0 & 0 \\
0 & \smash{\tfrac12} & 0 & \smash{\tfrac12} & 0 & \smash{\tfrac12} & \smash{\tfrac12} & 0 & \smash{\tfrac12} & \smash{\tfrac12} & 0 & 0 \\
0 & 1 & 0 & 1 & 0 & 0 & 0 & 0 & 0 & 0 & 0 & 1 \\
0 & 1 & 0 & 0 & 0 & 0 & 0 & 1 & 1 & 0 & 0 & 0 \\
0 & 0 & \smash{\tfrac12} & 0 & \smash{\tfrac12} & \smash{\tfrac12} & 0 & \smash{\tfrac12} & \smash{\tfrac12} & 0 & \smash{\tfrac12} & 0 \\
0 & 0 & 1 & 0 & 1 & 0 & 0 & 0 & 0 & 0 & 0 & 1 \\
0 & 0 & 1 & 0 & 0 & 0 & 1 & 0 & 1 & 0 & 0 & 0 \\
0 & 0 & 0 & 1 & 0 & 1 & 0 & 0 & 0 & 0 & 1 & 0 \\
0 & 0 & 0 & 0 & 1 & 1 & 0 & 0 & 0 & 1 & 0 & 0
\end{array}
\]
Eight of the twelve vertices of $P_0$ and eight of the nine vertices of
$P_1$ lie on the hyperplane $x_{e_1}=0$. These are the vertices of the
common facet along which the two polytopes are glued
(Example~\ref{ex:glue}), and the identification of
Theorem~\ref{thm:gluing} carries the first eight onto the second.
\end{example}

\section{Facets of the full-dimensional polytopes and Whitehead-move gluing}
\label{sec:facets}

Ten of the twenty types have $\dim P_i=6$: types
$0,1,2,4,6,7,12,13,17,19$. For these, $P_i$ is a $6$-dimensional
polytope in the affine space $\{A_ix=\mathbf 1\}$ whose facets lie in
the hyperplanes $x_e=0$. We ask whether its facets correspond to the
\WM{}s available at $G_i$, and whether two polytopes related by a
\WM{} meet along a common facet.

\subsection{Facets versus forced strata}

A hyperplane $x_e=0$ cuts $P_i$ in a \emph{facet} (a face of dimension
$5$) only if no further coordinate is forced on it; otherwise
$F_e=\{x\in P_i:x_e=0\}$ is a face of dimension at most $4$, whose
dimension we compute as in Proposition~\ref{prop:dimformula}, applied
to the system $A_ix=\mathbf1,\ x_e=0$.

\begin{example}
\label{ex:notfacet}
For type $1$, the row equation of the cone point $5$ is
$x_{e_6}+x_{e_7}+x_{e_8}+x_{e_{12}}=1$, and half of the combination
$(R_1+R_2+R_5+R_6)-(R_3+R_4)$ of the row equations is
\[
  x_{e_1}+x_{e_7}+x_{e_8}+x_{e_{10}}+x_{e_{11}}+x_{e_{12}}=1 .
\]
On $x_{e_6}=0$ the first equation gives $x_{e_7}+x_{e_8}+x_{e_{12}}=1$,
and the second then gives $x_{e_1}+x_{e_{10}}+x_{e_{11}}=0$, so
$x_{e_1}=x_{e_{10}}=x_{e_{11}}=0$. The face $F_{e_6}$ has dimension
$3$, not $5$, although the \WM{} at $e_6$ is admissible
(Table~\ref{tab:facetadm}).
\end{example}

\begin{table}[htbp]
\centering
\begin{tabular}{@{}c c c c c c@{}}
\toprule
type & \#facets & \#admissible edges & facet and admissible & facet only & admissible only \\
\midrule
0  & 12 & 12 & 12 & 0 & 0 \\
1  & 10 & 12 & 10 & 0 & 2 \\
2  & 7  & 10 & 7  & 0 & 3 \\
4  & 10 & 12 & 10 & 0 & 2 \\
6  & 7  & 10 & 7  & 0 & 3 \\
7  & 8  & 8  & 8  & 0 & 0 \\
12 & 8  & 8  & 8  & 0 & 0 \\
13 & 7  & 10 & 7  & 0 & 3 \\
17 & 8  & 8  & 8  & 0 & 0 \\
19 & 8  & 8  & 8  & 0 & 0 \\
\bottomrule
\end{tabular}
\caption{For each of the ten full-dimensional types: the number of
facets ($x_e=0$ cutting a $5$-dimensional face), the number of
edges at which the \WM{} is admissible, and their overlap. ``Facet
only'' is zero in every case (Theorem~\ref{thm:facetsubset});
``admissible only'' counts edges where the move is combinatorially
valid but $F_e$ has dimension less than $5$, as in
Example~\ref{ex:notfacet}.}
\label{tab:facetadm}
\end{table}

A \WM{} at the edge $e^*$ of $G$ is the diagonal flip at the dual edge
$e$ of $T$: delete $e$ from the quadrilateral formed by its two
triangles and insert the other diagonal. Following~\cite[Section~8.1]{Paper1}
it is \emph{admissible} when the resulting graph is again a CE-graph,
bigons allowed: connected and without a monogon. The flip lowers the
degrees of the two endpoints of $e$ by one each (by two for a loop),
raises the degrees of the two opposite corners and changes nothing
else, so it is admissible exactly when both endpoints of $e$ have
degree at least three, or, for a loop, its vertex has degree at least
four. Every loop of the twenty types is at a vertex of degree at
least seven, so only edges at cone points of degree two are
inadmissible.

\begin{lemma}
\label{lem:degreetwo}
Let $P_i$ be full-dimensional and let $e$ be an edge of $T_i$ with an
endpoint $v$ of degree two. Then $F_e$ is not a facet of $P_i$.
\end{lemma}

\begin{proof}
Let $f$ be the other edge at $v$, and write $e=va$, $f=vb$ (possibly
$a=b$); neither is a loop, since a loop at $v$ would be its only edge.
The two triangles at $v$ both have sides $e$ and $f$, so their third
sides $g,g'$ both join $a$ and $b$, and $g\ne g'$, since otherwise the
two triangles would close up into a sphere with three vertices. The
columns of $g$ and $g'$ in $A_i$ are equal, and $g,g'\notin\{e,f\}$.

On $F_e$ the row equation of $v$ gives $x_f=1$, and the row equation
of $b$ then forces every other coordinate at $b$ to vanish, including
$x_g$. If $F_e$ were a facet, its affine hull would be
$\aff(P_i)\cap\{x_e=0\}$. This affine space contains, with any of its
points $y$, the points $y+t(\mathbf 1_g-\mathbf 1_{g'})$ for all
$t\in\RR$, because $A_i(\mathbf 1_g-\mathbf 1_{g'})=0$ and the $e$-th
coordinate is unchanged. On these points $x_g$ is not constant,
contradicting $x_g=0$ on $F_e$.
\end{proof}

\begin{theorem}[Facets are Whitehead-admissible]
\label{thm:facetsubset}
For each of the ten full-dimensional types, every facet $x_e=0$ of
$P_i$ occurs at an edge $e$ at which the \WM{} is admissible. The
converse fails at $13$ of the $98$ admissible edges.
\end{theorem}

\begin{proof}
By the admissibility criterion above, an inadmissible edge has an
endpoint of degree two, and Lemma~\ref{lem:degreetwo} applies. The
counts are those of Table~\ref{tab:facetadm}, computed for all $120$
edges of the ten types; at the $13$ admissible edges that are not
facets, $F_e$ has dimension $1$, $2$ or $3$.
\end{proof}

\subsection{Whitehead moves glue matching facets}

\begin{theorem}[Gluing]
\label{thm:gluing}
Let $P_i$ be full-dimensional, let the \WM{} at $e$ be admissible, and
let it produce the type $j$, with $e'$ the new edge of $T_j$. Then
$T_i\smallsetminus e=T_j\smallsetminus e'$, dually the contracted
graphs $G_i/e^*$ and $G_j/e'^*$ are isomorphic, and dropping the
coordinates $e$ and $e'$ identifies the faces $F_e\subset P_i$ and
$F_{e'}\subset P_j$ with the same polytope, that of
$T_i\smallsetminus e$. In particular $F_e$ is a facet of $P_i$ if and
only if $F_{e'}$ is a facet of $P_j$.
\end{theorem}

\begin{proof}
The flip deletes $e$ and inserts $e'$, so the two triangulations agree
after deleting these edges; deleting an edge of $T$ is dual to
contracting the corresponding edge of $G$, the operation of the
contraction hierarchy~\cite[Definition~8.4]{Paper1}. On $F_e$ the
column of $e$ is multiplied by $0$, so dropping the coordinate $e$
maps $F_e$ bijectively onto
$\{y\ge0 : \text{incidence}(T_i\smallsetminus e)\,y=\mathbf 1\}$, and
likewise for $F_{e'}$. If $F_e$ is a facet, this polytope has
dimension $5$, so $F_{e'}$ has dimension $5$ as well, and it is a
facet of $P_j$ provided $P_j$ is full-dimensional; that it is, for
every facet, is Proposition~\ref{prop:closure} below.
\end{proof}

We also checked Theorem~\ref{thm:gluing} computationally, by
performing each of the $85$ flips at the dart level, identifying the
target type by a canonical form, and confirming that the facet of the
source and the face of the target agree, and that the inverse flip
returns the source.

\begin{example}[The cube and the stacked type glue along $e_1$]
\label{ex:glue}
Contracting $G_0$ along $e_1$ gives the seven-vertex graph $H$ with
degree sequence $(3,3,3,3,3,3,4)$, face signature $(3,3,4,4,4,4)$ and
edge list
\begin{align*}
e_1&{=}(1,2), & e_2&{=}(1,3), & e_3&{=}(1,4), & e_4&{=}(1,5),\\
e_5&{=}(2,4), & e_6&{=}(2,6), & e_7&{=}(3,5), & e_8&{=}(3,6),\\
e_9&{=}(4,7), & e_{10}&{=}(5,7), & e_{11}&{=}(6,7). &&
\end{align*}
The \WM{} at $e_1$ sends $G_0$ to $G_1$, and the new edge is the edge
$e_1$ of $G_1$ in the labeling of Table~\ref{tab:edgelists};
contracting $G_1$ along it gives the same $H$. The facets
$\{x\in P_0:x_{e_1}=0\}$ and $\{x\in P_1:x_{e_1}=0\}$ are both, after
dropping the zero coordinate, the polytope of $H^*$, and $P_0$ and
$P_1$ are glued along it (compare
Examples~\ref{ex:cubevertices}--\ref{ex:stackedvertices}).
\end{example}

\subsection{The facet-gluing graph}

\begin{proposition}
\label{prop:closure}
For every facet of every one of the ten full-dimensional types, the
\WM{} at that edge leads to one of the same ten types; none leads to
a lower-dimensional or infeasible type.
\end{proposition}

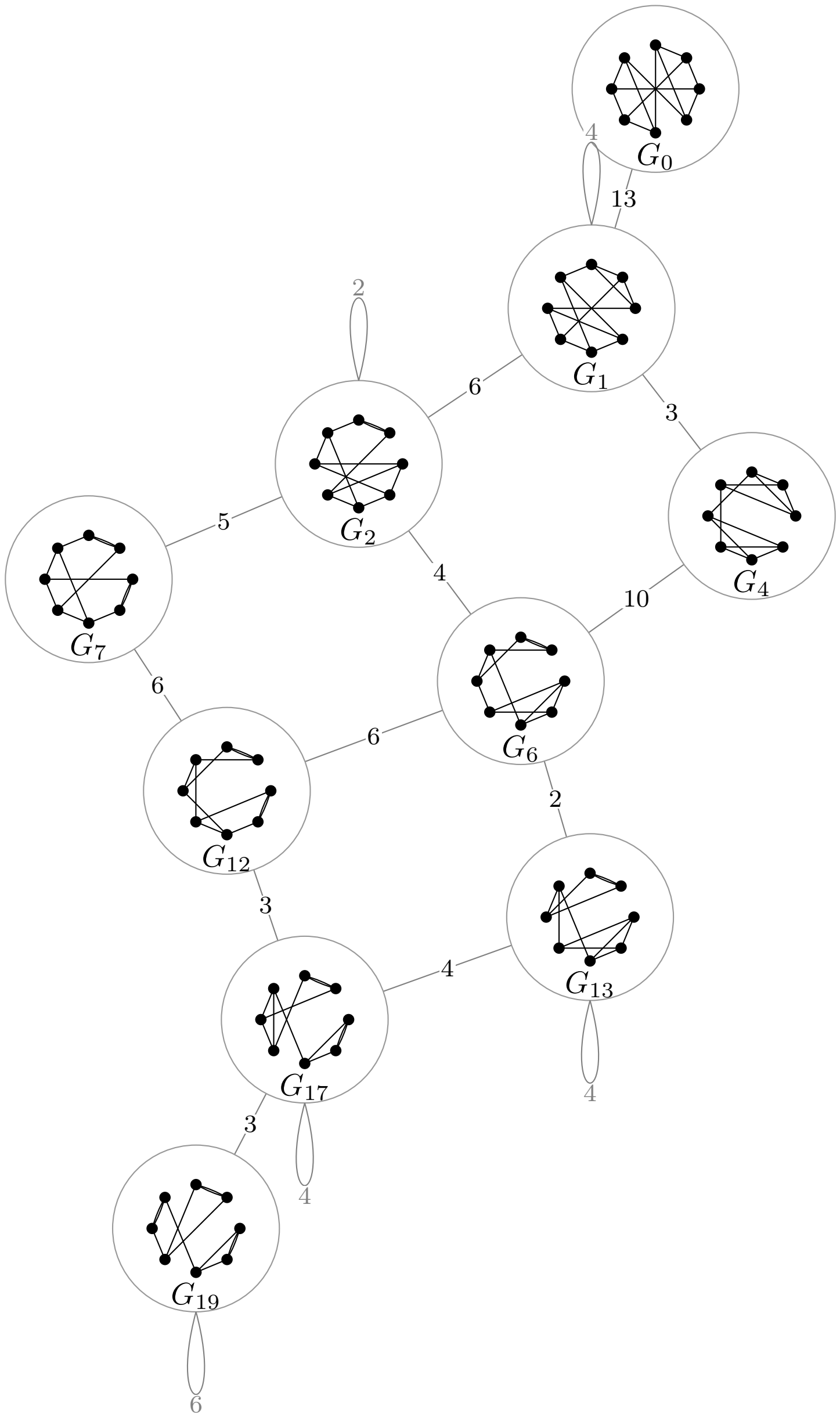
\begin{figure}[p]
\centering
\begin{tikzpicture}[vtx/.style={circle,fill=black,inner sep=1.3pt}]
  \draw[gray] (2.587,6.500) -- (1.858,4.000);
  \node[fill=white,inner sep=0.5pt,font=\scriptsize] at (2.222,5.250) {13};
  \draw[gray] (1.858,4.000) -- (-0.795,2.227);
  \node[fill=white,inner sep=0.5pt,font=\scriptsize] at (0.531,3.114) {6};
  \draw[gray] (1.858,4.000) -- (3.684,1.634);
  \node[fill=white,inner sep=0.5pt,font=\scriptsize] at (2.771,2.817) {3};
  \draw[gray] (-0.795,2.227) -- (1.052,-0.248);
  \node[fill=white,inner sep=0.5pt,font=\scriptsize] at (0.128,0.990) {4};
  \draw[gray] (-0.795,2.227) -- (-3.871,0.913);
  \node[fill=white,inner sep=0.5pt,font=\scriptsize] at (-2.333,1.570) {5};
  \draw[gray] (3.684,1.634) -- (1.052,-0.248);
  \node[fill=white,inner sep=0.5pt,font=\scriptsize] at (2.368,0.693) {10};
  \draw[gray] (1.052,-0.248) -- (-2.297,-1.498);
  \node[fill=white,inner sep=0.5pt,font=\scriptsize] at (-0.623,-0.873) {6};
  \draw[gray] (1.052,-0.248) -- (1.841,-2.937);
  \node[fill=white,inner sep=0.5pt,font=\scriptsize] at (1.446,-1.593) {2};
  \draw[gray] (-3.871,0.913) -- (-2.297,-1.498);
  \node[fill=white,inner sep=0.5pt,font=\scriptsize] at (-3.084,-0.293) {6};
  \draw[gray] (-2.297,-1.498) -- (-1.410,-4.105);
  \node[fill=white,inner sep=0.5pt,font=\scriptsize] at (-1.854,-2.802) {3};
  \draw[gray] (1.841,-2.937) -- (-1.410,-4.105);
  \node[fill=white,inner sep=0.5pt,font=\scriptsize] at (0.216,-3.521) {4};
  \draw[gray] (-1.410,-4.105) -- (-2.649,-6.485);
  \node[fill=white,inner sep=0.5pt,font=\scriptsize] at (-2.029,-5.295) {3};
  \begin{scope}[shift={(2.587,6.500)}]
    \fill[white] (0,0) circle (0.95);
    \draw[black!40] (0,0) circle (0.95);
    \begin{scope}[scale=0.5]
  \draw (0.000,1.000) -- (0.707,0.707);
  \draw (0.707,0.707) -- (1.000,0.000);
  \draw (0.707,-0.707) -- (1.000,0.000);
  \draw (0.000,1.000) -- (0.707,-0.707);
  \draw (0.000,-1.000) -- (-0.707,-0.707);
  \draw (-0.707,-0.707) -- (-1.000,-0.000);
  \draw (-0.707,0.707) -- (-1.000,-0.000);
  \draw (0.000,-1.000) -- (-0.707,0.707);
  \draw (0.000,1.000) -- (0.000,-1.000);
  \draw (0.707,0.707) -- (-0.707,-0.707);
  \draw (1.000,0.000) -- (-1.000,-0.000);
  \draw (0.707,-0.707) -- (-0.707,0.707);
  \node[vtx] at (0.000,1.000) {};
  \node[vtx] at (0.707,0.707) {};
  \node[vtx] at (1.000,0.000) {};
  \node[vtx] at (0.707,-0.707) {};
  \node[vtx] at (0.000,-1.000) {};
  \node[vtx] at (-0.707,-0.707) {};
  \node[vtx] at (-1.000,-0.000) {};
  \node[vtx] at (-0.707,0.707) {};
    \end{scope}
    \node[font=\bfseries\small] at (0,-0.77) {$G_{0}$};
  \end{scope}
  \begin{scope}[shift={(1.858,4.000)}]
    \fill[white] (0,0) circle (0.95);
    \draw[black!40] (0,0) circle (0.95);
    \begin{scope}[scale=0.5]
  \draw (0.000,1.000) -- (0.707,0.707);
  \draw (0.707,0.707) -- (1.000,0.000);
  \draw (0.000,1.000) -- (1.000,0.000);
  \draw (0.707,-0.707) -- (0.000,-1.000);
  \draw (0.000,-1.000) -- (-0.707,-0.707);
  \draw (-0.707,-0.707) -- (-1.000,-0.000);
  \draw (0.707,-0.707) -- (-1.000,-0.000);
  \draw (-0.707,0.707) -- (0.707,-0.707);
  \draw (-0.707,0.707) -- (0.000,-1.000);
  \draw (0.707,0.707) -- (-0.707,-0.707);
  \draw (1.000,0.000) -- (-1.000,-0.000);
  \draw (0.000,1.000) -- (-0.707,0.707);
  \node[vtx] at (0.000,1.000) {};
  \node[vtx] at (0.707,0.707) {};
  \node[vtx] at (1.000,0.000) {};
  \node[vtx] at (0.707,-0.707) {};
  \node[vtx] at (0.000,-1.000) {};
  \node[vtx] at (-0.707,-0.707) {};
  \node[vtx] at (-1.000,-0.000) {};
  \node[vtx] at (-0.707,0.707) {};
    \end{scope}
    \node[font=\bfseries\small] at (0,-0.77) {$G_{1}$};
  \end{scope}
  \begin{scope}[shift={(-0.795,2.227)}]
    \fill[white] (0,0) circle (0.95);
    \draw[black!40] (0,0) circle (0.95);
    \begin{scope}[scale=0.5]
  \draw (0.000,1.000) to[bend left=-6] (0.707,0.707);
  \draw (0.000,1.000) to[bend left=6] (0.707,0.707);
  \draw (1.000,0.000) -- (0.707,-0.707);
  \draw (0.000,-1.000) -- (0.707,-0.707);
  \draw (0.000,-1.000) -- (-0.707,-0.707);
  \draw (-0.707,-0.707) -- (1.000,0.000);
  \draw (-1.000,-0.000) -- (1.000,0.000);
  \draw (-1.000,-0.000) -- (0.707,-0.707);
  \draw (-0.707,0.707) -- (0.000,-1.000);
  \draw (0.707,0.707) -- (-0.707,-0.707);
  \draw (-0.707,0.707) -- (-1.000,-0.000);
  \draw (0.000,1.000) -- (-0.707,0.707);
  \node[vtx] at (0.000,1.000) {};
  \node[vtx] at (0.707,0.707) {};
  \node[vtx] at (1.000,0.000) {};
  \node[vtx] at (0.707,-0.707) {};
  \node[vtx] at (0.000,-1.000) {};
  \node[vtx] at (-0.707,-0.707) {};
  \node[vtx] at (-1.000,-0.000) {};
  \node[vtx] at (-0.707,0.707) {};
    \end{scope}
    \node[font=\bfseries\small] at (0,-0.77) {$G_{2}$};
  \end{scope}
  \begin{scope}[shift={(3.684,1.634)}]
    \fill[white] (0,0) circle (0.95);
    \draw[black!40] (0,0) circle (0.95);
    \begin{scope}[scale=0.5]
  \draw (0.000,1.000) -- (0.707,0.707);
  \draw (0.707,0.707) -- (1.000,0.000);
  \draw (0.000,1.000) -- (1.000,0.000);
  \draw (0.707,-0.707) -- (0.000,-1.000);
  \draw (0.000,-1.000) -- (-0.707,-0.707);
  \draw (0.707,-0.707) -- (-0.707,-0.707);
  \draw (-1.000,-0.000) -- (0.707,-0.707);
  \draw (-1.000,-0.000) -- (0.000,-1.000);
  \draw (0.707,0.707) -- (-0.707,0.707);
  \draw (1.000,0.000) -- (-0.707,0.707);
  \draw (0.000,1.000) -- (-1.000,-0.000);
  \draw (-0.707,0.707) -- (-0.707,-0.707);
  \node[vtx] at (0.000,1.000) {};
  \node[vtx] at (0.707,0.707) {};
  \node[vtx] at (1.000,0.000) {};
  \node[vtx] at (0.707,-0.707) {};
  \node[vtx] at (0.000,-1.000) {};
  \node[vtx] at (-0.707,-0.707) {};
  \node[vtx] at (-1.000,-0.000) {};
  \node[vtx] at (-0.707,0.707) {};
    \end{scope}
    \node[font=\bfseries\small] at (0,-0.77) {$G_{4}$};
  \end{scope}
  \begin{scope}[shift={(1.052,-0.248)}]
    \fill[white] (0,0) circle (0.95);
    \draw[black!40] (0,0) circle (0.95);
    \begin{scope}[scale=0.5]
  \draw (0.000,1.000) to[bend left=-6] (0.707,0.707);
  \draw (0.000,1.000) to[bend left=6] (0.707,0.707);
  \draw (1.000,0.000) -- (0.707,-0.707);
  \draw (0.707,-0.707) -- (0.000,-1.000);
  \draw (1.000,0.000) -- (0.000,-1.000);
  \draw (-0.707,-0.707) -- (1.000,0.000);
  \draw (-0.707,-0.707) -- (0.707,-0.707);
  \draw (-1.000,-0.000) -- (-0.707,0.707);
  \draw (0.707,0.707) -- (-0.707,0.707);
  \draw (-1.000,-0.000) -- (-0.707,-0.707);
  \draw (0.000,1.000) -- (-1.000,-0.000);
  \draw (-0.707,0.707) -- (0.000,-1.000);
  \node[vtx] at (0.000,1.000) {};
  \node[vtx] at (0.707,0.707) {};
  \node[vtx] at (1.000,0.000) {};
  \node[vtx] at (0.707,-0.707) {};
  \node[vtx] at (0.000,-1.000) {};
  \node[vtx] at (-0.707,-0.707) {};
  \node[vtx] at (-1.000,-0.000) {};
  \node[vtx] at (-0.707,0.707) {};
    \end{scope}
    \node[font=\bfseries\small] at (0,-0.77) {$G_{6}$};
  \end{scope}
  \begin{scope}[shift={(-3.871,0.913)}]
    \fill[white] (0,0) circle (0.95);
    \draw[black!40] (0,0) circle (0.95);
    \begin{scope}[scale=0.5]
  \draw (0.000,1.000) to[bend left=-6] (0.707,0.707);
  \draw (0.000,1.000) to[bend left=6] (0.707,0.707);
  \draw (1.000,0.000) to[bend left=-6] (0.707,-0.707);
  \draw (1.000,0.000) to[bend left=6] (0.707,-0.707);
  \draw (0.000,-1.000) -- (0.707,-0.707);
  \draw (0.000,-1.000) -- (-0.707,-0.707);
  \draw (-1.000,-0.000) -- (-0.707,-0.707);
  \draw (-0.707,0.707) -- (0.000,-1.000);
  \draw (0.707,0.707) -- (-0.707,-0.707);
  \draw (-0.707,0.707) -- (-1.000,-0.000);
  \draw (0.000,1.000) -- (-0.707,0.707);
  \draw (-1.000,-0.000) -- (1.000,0.000);
  \node[vtx] at (0.000,1.000) {};
  \node[vtx] at (0.707,0.707) {};
  \node[vtx] at (1.000,0.000) {};
  \node[vtx] at (0.707,-0.707) {};
  \node[vtx] at (0.000,-1.000) {};
  \node[vtx] at (-0.707,-0.707) {};
  \node[vtx] at (-1.000,-0.000) {};
  \node[vtx] at (-0.707,0.707) {};
    \end{scope}
    \node[font=\bfseries\small] at (0,-0.77) {$G_{7}$};
  \end{scope}
  \begin{scope}[shift={(-2.297,-1.498)}]
    \fill[white] (0,0) circle (0.95);
    \draw[black!40] (0,0) circle (0.95);
    \begin{scope}[scale=0.5]
  \draw (0.000,1.000) to[bend left=-6] (0.707,0.707);
  \draw (0.000,1.000) to[bend left=6] (0.707,0.707);
  \draw (0.707,-0.707) to[bend left=-6] (1.000,0.000);
  \draw (0.707,-0.707) to[bend left=6] (1.000,0.000);
  \draw (0.000,-1.000) -- (-0.707,-0.707);
  \draw (0.000,-1.000) -- (0.707,-0.707);
  \draw (-1.000,-0.000) -- (-0.707,0.707);
  \draw (0.707,0.707) -- (-0.707,0.707);
  \draw (-1.000,-0.000) -- (0.000,-1.000);
  \draw (0.000,1.000) -- (-1.000,-0.000);
  \draw (-0.707,0.707) -- (-0.707,-0.707);
  \draw (-0.707,-0.707) -- (1.000,0.000);
  \node[vtx] at (0.000,1.000) {};
  \node[vtx] at (0.707,0.707) {};
  \node[vtx] at (1.000,0.000) {};
  \node[vtx] at (0.707,-0.707) {};
  \node[vtx] at (0.000,-1.000) {};
  \node[vtx] at (-0.707,-0.707) {};
  \node[vtx] at (-1.000,-0.000) {};
  \node[vtx] at (-0.707,0.707) {};
    \end{scope}
    \node[font=\bfseries\small] at (0,-0.77) {$G_{12}$};
  \end{scope}
  \begin{scope}[shift={(1.841,-2.937)}]
    \fill[white] (0,0) circle (0.95);
    \draw[black!40] (0,0) circle (0.95);
    \begin{scope}[scale=0.5]
  \draw (0.000,1.000) to[bend left=-6] (0.707,0.707);
  \draw (0.000,1.000) to[bend left=6] (0.707,0.707);
  \draw (1.000,0.000) -- (0.707,-0.707);
  \draw (0.000,-1.000) -- (0.707,-0.707);
  \draw (0.000,-1.000) -- (1.000,0.000);
  \draw (-0.707,-0.707) -- (1.000,0.000);
  \draw (-0.707,-0.707) -- (0.707,-0.707);
  \draw (-1.000,-0.000) -- (0.707,0.707);
  \draw (-0.707,0.707) -- (-0.707,-0.707);
  \draw (0.000,1.000) -- (-1.000,-0.000);
  \draw (-0.707,0.707) -- (0.000,-1.000);
  \draw (-1.000,-0.000) -- (-0.707,0.707);
  \node[vtx] at (0.000,1.000) {};
  \node[vtx] at (0.707,0.707) {};
  \node[vtx] at (1.000,0.000) {};
  \node[vtx] at (0.707,-0.707) {};
  \node[vtx] at (0.000,-1.000) {};
  \node[vtx] at (-0.707,-0.707) {};
  \node[vtx] at (-1.000,-0.000) {};
  \node[vtx] at (-0.707,0.707) {};
    \end{scope}
    \node[font=\bfseries\small] at (0,-0.77) {$G_{13}$};
  \end{scope}
  \begin{scope}[shift={(-1.410,-4.105)}]
    \fill[white] (0,0) circle (0.95);
    \draw[black!40] (0,0) circle (0.95);
    \begin{scope}[scale=0.5]
  \draw (0.000,1.000) to[bend left=-6] (0.707,0.707);
  \draw (0.000,1.000) to[bend left=6] (0.707,0.707);
  \draw (1.000,0.000) to[bend left=-6] (0.707,-0.707);
  \draw (1.000,0.000) to[bend left=6] (0.707,-0.707);
  \draw (0.000,-1.000) -- (0.707,-0.707);
  \draw (-0.707,-0.707) -- (-1.000,-0.000);
  \draw (0.707,0.707) -- (-1.000,-0.000);
  \draw (-0.707,-0.707) -- (-0.707,0.707);
  \draw (0.000,1.000) -- (-0.707,-0.707);
  \draw (-1.000,-0.000) -- (-0.707,0.707);
  \draw (0.000,-1.000) -- (1.000,0.000);
  \draw (-0.707,0.707) -- (0.000,-1.000);
  \node[vtx] at (0.000,1.000) {};
  \node[vtx] at (0.707,0.707) {};
  \node[vtx] at (1.000,0.000) {};
  \node[vtx] at (0.707,-0.707) {};
  \node[vtx] at (0.000,-1.000) {};
  \node[vtx] at (-0.707,-0.707) {};
  \node[vtx] at (-1.000,-0.000) {};
  \node[vtx] at (-0.707,0.707) {};
    \end{scope}
    \node[font=\bfseries\small] at (0,-0.77) {$G_{17}$};
  \end{scope}
  \begin{scope}[shift={(-2.649,-6.485)}]
    \fill[white] (0,0) circle (0.95);
    \draw[black!40] (0,0) circle (0.95);
    \begin{scope}[scale=0.5]
  \draw (0.000,1.000) to[bend left=-6] (0.707,0.707);
  \draw (0.000,1.000) to[bend left=6] (0.707,0.707);
  \draw (1.000,0.000) to[bend left=-6] (0.707,-0.707);
  \draw (1.000,0.000) to[bend left=6] (0.707,-0.707);
  \draw (0.000,-1.000) -- (0.707,-0.707);
  \draw (-0.707,-0.707) -- (0.707,0.707);
  \draw (-1.000,-0.000) to[bend left=-6] (-0.707,0.707);
  \draw (-1.000,-0.000) to[bend left=6] (-0.707,0.707);
  \draw (0.000,1.000) -- (-0.707,-0.707);
  \draw (0.000,-1.000) -- (1.000,0.000);
  \draw (-0.707,0.707) -- (0.000,-1.000);
  \draw (-0.707,-0.707) -- (-1.000,-0.000);
  \node[vtx] at (0.000,1.000) {};
  \node[vtx] at (0.707,0.707) {};
  \node[vtx] at (1.000,0.000) {};
  \node[vtx] at (0.707,-0.707) {};
  \node[vtx] at (0.000,-1.000) {};
  \node[vtx] at (-0.707,-0.707) {};
  \node[vtx] at (-1.000,-0.000) {};
  \node[vtx] at (-0.707,0.707) {};
    \end{scope}
    \node[font=\bfseries\small] at (0,-0.77) {$G_{19}$};
  \end{scope}
  \draw[gray] (1.858,4.950) to[loop above,looseness=10,min distance=1.3cm] node[font=\scriptsize,fill=white,inner sep=0.5pt] {4} (1.858,4.950);
  \draw[gray] (-0.795,3.177) to[loop above,looseness=10,min distance=1.3cm] node[font=\scriptsize,fill=white,inner sep=0.5pt] {2} (-0.795,3.177);
  \draw[gray] (1.841,-3.887) to[loop below,looseness=10,min distance=1.3cm] node[font=\scriptsize,fill=white,inner sep=0.5pt] {4} (1.841,-3.887);
  \draw[gray] (-1.410,-5.055) to[loop below,looseness=10,min distance=1.3cm] node[font=\scriptsize,fill=white,inner sep=0.5pt] {4} (-1.410,-5.055);
  \draw[gray] (-2.649,-7.435) to[loop below,looseness=10,min distance=1.3cm] node[font=\scriptsize,fill=white,inner sep=0.5pt] {6} (-2.649,-7.435);
\end{tikzpicture}
\caption{The facet-gluing graph of Theorem~\ref{thm:gluing}. One node
per full-dimensional type, drawn as its Voronoi graph $G_i$; an edge
joins two types when a facet of one is glued to a facet of the other,
labeled by the number of glued facets counted from both sides (so the
$13$ between $G_0$ and $G_1$ counts the twelve facets of $P_0$ and the
one facet of $P_1$ that lead to each other), and a loop at $G_i$
counts the facets of $P_i$ glued to $P_i$ itself. The labels sum to
$85$, the number of facets. Node positions are a Kamada--Kawai layout
and carry no meaning.}
\label{fig:gluingmeta}
\end{figure}

\begin{proof}[Proof (computational)]
All $85$ targets lie in $\{0,1,2,4,6,7,12,13,17,19\}$
(Figure~\ref{fig:gluingmeta}).
\end{proof}

By Theorems~\ref{thm:facetsubset} and~\ref{thm:gluing} and
Proposition~\ref{prop:closure}, every facet of every full-dimensional
cell is identified with a facet of a full-dimensional cell, and the
identifications are recorded by the graph of
Figure~\ref{fig:gluingmeta}. This is the top-dimensional part of the
structure suggested in~\cite[Section~8]{Paper1} by the contraction
hierarchy, in which the walls between top cells are the level-one
graphs. The ``admissible only'' edges of Table~\ref{tab:facetadm} show
that this gluing is a proper part of the full \WM{} orbit
of~\cite[Theorem~8.1]{Paper1}: combinatorial adjacency in the orbit
does not always come with a common wall, because degenerating the
edge $e$ forces other edges to degenerate with it
(Example~\ref{ex:notfacet}). Section~\ref{sec:unfolding} shows that
the identifications are compatible with the automorphism groups of
the cells.

\subsection{Generating the vertices from the cube}
\label{sec:vertexgen}

The identification of Theorem~\ref{thm:gluing} carries vertices
across, and iterating it over the graph of Figure~\ref{fig:gluingmeta}
generates almost every vertex of every full-dimensional type from the
twelve vertices of $P_0$ (Example~\ref{ex:cubevertices}).

\begin{definition}[Vertex transfer]
\label{def:transfer}
Let the \WM{} at $e$ send $T_i$ to $T_j$ with new edge $e'$, and let
$\sigma$ be the resulting bijection from the edges of $T_i$ other
than $e$ to the edges of $T_j$ other than $e'$, composed with the
relabeling of $T_j$ into the labeling of Table~\ref{tab:edgelists}
(defined up to $\Aut^+(T_j)$). For $x\in P_i$ with $x_e=0$ define
$\tau(x)_{e'}=0$ and $\tau(x)_{\sigma(f)}=x_f$ for every other edge
$f$.
\end{definition}

\begin{proposition}[Transfer preserves vertices]
\label{prop:vertextransfer}
If $x$ is a vertex of $P_i$ with $x_e=0$, then $\tau(x)$ is a vertex of
$P_j$.
\end{proposition}

\begin{proof}
By Theorem~\ref{thm:gluing}, $\tau$ restricts to an affine
isomorphism from the face $F_e$ of $P_i$ onto the face $F_{e'}$ of
$P_j$. It maps vertices of $F_e$ to vertices of $F_{e'}$, and a vertex
of a face of a polytope is a vertex of the polytope.
\end{proof}

\begin{theorem}[Vertex generation from the cube]
\label{thm:vertexgen}
Starting from the twelve vertices of $P_0$ and applying all transfers
along the $85$ facet identifications, in either direction and for
every choice of $\sigma$, until nothing new appears, generates every
vertex of $P_0,P_1,P_2,P_4,P_6,P_7,P_{12}$ (the seven full-dimensional
types whose $T$ has no loop), all but one of the seven vertices of
$P_{13}$, all but two of the nine of $P_{17}$, and all but five of
the eleven of $P_{19}$. The eight vertices never generated are exactly
the vertices at which some loop of $T$ has the value $\tfrac12$.
Allowing transfers through the lower-dimensional types as well does
not change this.
\end{theorem}

\begin{proof}
The positive statements are a finite computation, which we carried
out. For the negative one, a flip changes the endpoints of no edge
other than $e$, so a loop of $T_j$ other than the new edge $e'$ is,
under $\sigma$, a loop of $T_i$ with the same coordinate, while
$\tau(x)_{e'}=0$. The cube has no loops, so by induction every
generated vertex is $0$ on every loop.
\end{proof}

So the transfer can only carry structure present in the seed; the
vertices of the triangulations with loops that use a loop are the
only ones it cannot reach.

\section{Automorphisms and the free unfolding}
\label{sec:unfolding}

\subsection{The gluing is compatible with the automorphism groups}

Each of the ten full-dimensional types appears in
Figure~\ref{fig:gluingmeta} as a single polytope, but a type can have
several facets glued to the same neighbour: all twelve facets of the
cube $P_0$ lead to type $1$, while only one facet of $P_1$ leads back
to type $0$. This is accounted for by symmetry. The group
$\Aut^+(T_i)$ of orientation-preserving automorphisms of the
triangulation acts on $P_i$ by permuting coordinates and permutes its
facets; for the cube it acts transitively on all twelve.

\begin{proposition}
\label{prop:orbitpairing}
For each of the $85$ facets, the identification of
Theorem~\ref{thm:gluing} maps the $\Aut^+(T_i)$-orbit of the facet
onto a single $\Aut^+(T_j)$-orbit of facets of the target, and the
resulting map on orbits is a fixed-point-free involution. The ten
cells have $1,6,7,3,7,3,5,7,8,5$ facet orbits (for the types
$0,1,2,4,6,7,12,13,17,19$), $52$ in all, which are paired into $26$
walls.
\end{proposition}

\begin{proof}[Proof (computational)]
We computed $\Aut^+(T_i)$ as the set of dart permutations realizing the
minimal canonical code, of orders $24,2,1,4,1,4,2,1,1,2$, and checked
the statement directly.
\end{proof}

Consequently the quotients $P_i/\Aut^+(T_i)$ of the ten
full-dimensional cells, glued along the $26$ walls, form a finite
complex with orbifold cells. This is the combinatorial model of the
top-dimensional part of the Delaunay decomposition of the moduli
space of six unordered cone points, and hence of genus-two surfaces
via the hyperelliptic correspondence~\cite[Proposition~5.2]{Paper1}.
Whether the wedge-angle coordinates identify each top cell of that
decomposition with $P_i^\circ/\Aut^+(T_i)$ is the realizability
question left open in~\cite[Section~7]{Paper1}.

\subsection{The free unfolding}
\label{sec:unfoldnevercloses}

Passing from the quotient to Teichm\"uller space requires one copy of
each cell for each marked Delaunay triangulation. A naive way to
build such copies is the \emph{free unfolding}: start from one copy
of $P_0$ and attach a fresh copy of the target cell to every open
facet, recursively, never identifying two copies. Its combinatorics
is a tree on \emph{arcs}: an arc is a pair $(t,e)$ of a type and one
of its $85$ facets, and a copy of type $t$ entered through the facet
$e$ has one child for each other facet of $t$.

\begin{table}[htbp]
\centering
\begin{tabular}{@{}c r r c@{}}
\toprule
depth & new copies & cumulative total & ratio to previous depth\\
\midrule
0 & 1         & 1         & --- \\
1 & 12        & 13        & 12.00 \\
2 & 108       & 121       & 9.00 \\
3 & 828       & 949       & 7.67 \\
4 & 6012      & 6961      & 7.26 \\
5 & 44292     & 51253     & 7.37 \\
6 & 322116    & 373369    & 7.27 \\
7 & 2327988   & 2701357   & 7.23 \\
\bottomrule
\end{tabular}
\caption{Growth of the free unfolding, starting from one copy of the
cube. All counts from depth one on are divisible by $12$, as they must
be, since $\Aut^+(T_0)$ permutes the twelve facets of the root
transitively.}
\label{tab:growth}
\end{table}

\begin{proposition}[Exponential growth]
\label{prop:growth}
The number of copies at depth $d$ in the free unfolding grows like
$\lambda^d$, where $\lambda\approx7.0928$ is the spectral radius of
the $85\times85$ arc transition matrix; the next eigenvalue has
modulus $\approx6.0438$.
\end{proposition}

The ratios of Table~\ref{tab:growth} approach $\lambda$ slowly; the
ratio at depth $40$ is $7.0934$.

\subsection{Relation to Teichm\"uller space}

The free unfolding is not a cell decomposition of Teichm\"uller
space, and its growth rate is not that of the mapping class group. In
the Delaunay decomposition of $\TT_{0,6}\cong\TT_2$, an instance of
the canonical decompositions of~\cite{BowditchEpstein1988,Penner1987},
the cells are indexed by marked Delaunay triangulations, and the cells around a
face of codimension two close up: if the two degenerate edges lie in
disjoint quadrilaterals there are four cells around it, and if they
lie in a common pentagon there are five. The free unfolding never
identifies copies, so its dual graph is a tree; each codimension-two
face lies in only two cells, and the unfolding is not a manifold. The
dual graph of the true decomposition is the realizable part of the
flip graph of the six-punctured sphere, on which the mapping class
group acts with finite quotient, namely the graph of
Figure~\ref{fig:gluingmeta} modulo the automorphisms. It too grows
exponentially, since the mapping class group of the six-punctured
sphere contains free subgroups, but at a rate we have not computed.

What the computation does show is consistent with the classical
picture: $\TT_2$ is homeomorphic to $\RR^6$, for instance by
Fenchel--Nielsen coordinates~\cite{FarbMargalit2012}, and a locally
finite cell decomposition of it by $6$-dimensional cells must be
infinite, with every wall shared by exactly two cells.
Theorems~\ref{thm:facetsubset} and~\ref{thm:gluing} and
Propositions~\ref{prop:closure} and~\ref{prop:orbitpairing} are the
finite combinatorial shadow of that structure.

\section{Computational verification and formal proofs}
\label{sec:verification}

\subsection{Independent computation}

The script \texttt{paper2\_verification.py} recomputes every finite
statement of this paper from scratch, without using the tables: it
enumerates the triangulations of the sphere with six vertices and no
vertex of degree one as dart systems (rotation systems), recovering
the twenty trivalent types and their signatures, and matches them to
the reference labeling of Table~\ref{tab:edgelists} by an explicit
isomorphism. It then checks the ranks
(Theorem~\ref{thm:rankalways6}); the forced coordinates, dimensions and
infeasibility (Table~\ref{tab:results}), by linear programming and
exact rational arithmetic; the row combinations of
Sections~\ref{subsec:example3}--\ref{subsec:example14} and
Example~\ref{ex:notfacet}; the structure of
$T_i\smallsetminus Z_i$ in Proposition~\ref{prop:degenerate}; the
vertices, their half-integrality and the subgraph count of
Corollary~\ref{cor:vertexdescription} (Table~\ref{tab:vertices} and
Examples~\ref{ex:cubevertices}--\ref{ex:stackedvertices}); the
facets, admissible edges, flip targets and the admissibility criterion
(Table~\ref{tab:facetadm} and Figure~\ref{fig:gluingmeta}); the
contraction of Example~\ref{ex:glue}; vertex generation
(Theorem~\ref{thm:vertexgen}); the orbit pairing
(Proposition~\ref{prop:orbitpairing}); and the growth of the free
unfolding (Table~\ref{tab:growth} and Proposition~\ref{prop:growth}).
It uses Python~3 with \texttt{numpy}, \texttt{scipy} and
\texttt{networkx}, runs in under a minute, and ends with
\texttt{ALL CHECKS PASSED}.

\subsection{Lean formalization}

The file \texttt{Paper2CirclePatterns.lean} formalizes the finite
linear algebra in Lean~4~\cite{Lean4} (version 4.15.0, core library
only). From the cone-point lists of Table~\ref{tab:edgelists} it
builds the twenty incidence matrices and proves:
\begin{enumerate}[label=(\roman*)]
\item $\rank A_i=6$ for all twenty types
(Theorem~\ref{thm:rankalways6});
\item that $A_{14}x=N\mathbf 1$ with $x\in\mathbb Z_{\ge0}^{12}$ forces
$N=0$, so $P_{14}$ has no rational point and is empty
(Section~\ref{subsec:example14});
\item for each of the nine lower-dimensional types, that every
$x\in\mathbb Z_{\ge0}^{12}$ with $A_ix=N\mathbf 1$ has the forced
values of Table~\ref{tab:results} (times $N$), which by clearing
denominators is the statement for rational points;
\item the dimensions of Theorem~\ref{thm:dimtable}, as an upper bound
from the rank of $[A_i;(e_j)_{j\in F}]$ and a matching lower bound from
the affine rank of the vertices;
\item that every basic feasible solution of every type is
half-integral and that the vertices are exactly those of
Table~\ref{tab:vertices}, with the stated counts and
$(m,h)$-distributions (Theorem~\ref{thm:halfintegral});
\item the facet counts and admissible-edge counts of
Table~\ref{tab:facetadm}, and that every facet is at an admissible
edge (Theorem~\ref{thm:facetsubset}).
\end{enumerate}
The forced coordinates and the infeasibility are proved by the
decision procedure \texttt{omega} and checked by the Lean kernel; the
other statements are finite computations checked with
\texttt{native\_decide}, which also trusts the Lean compiler. The
formalization uses, without proving them, the standard facts that the
vertices of $P_i$ are its basic feasible solutions and that $P_i$ is
their convex hull. The general statements, Lemma~\ref{lem:rankwithloops},
Theorem~\ref{thm:halfintegral} for arbitrary graphs and
Lemma~\ref{lem:degreetwo}, are proved in the text; formalizing them
would need a library of graph theory and polyhedral geometry such as
Mathlib. The flip computations (Theorem~\ref{thm:gluing},
Propositions~\ref{prop:closure} and~\ref{prop:orbitpairing}, and the
unfolding) are checked by the Python script only.

\section{Discussion}
\label{sec:discussion}

\subsection{Genericity does not predict dimension}

The full-dimensional types do not coincide with the generic ones
(simple $T$, types $0$ and $1$ only): eight non-generic types
($2,4,6,7,12,13,17,19$) are also full-dimensional, and among the six
types with loops, three ($13,17,19$) are full-dimensional and three
($15,16,18$) are not. By Proposition~\ref{prop:degenerate}(a) the
relevant dichotomy is realizability itself. Whether it can be read
off from a simple invariant of $T$, without solving the system, is
open. A necessary condition is visible in
Lemma~\ref{lem:degreetwo} and Section~\ref{subsec:example14}:
degree-two cone points force equalities, and four of them, as in type
$14$, leave no room at all.

\subsection{The meaning of a forced coordinate}

A forced coordinate $x_e=0$ means that every Delaunay triangulation
of the type has the four cone points around $e$ concyclic, so that
the Voronoi edge $e^*$ has length zero and $e$ is not an edge of the
Delaunay decomposition~\cite[Proposition~7.7]{Paper1}. A forced
coordinate $x_e=1$ means that a single wedge fills the whole cone
angle at an endpoint of $e$, so that the Voronoi cell of that cone
point would be a monogon, which~\cite[Theorem~6.5]{Paper1} excludes.
Proposition~\ref{prop:degenerate} turns these two observations into
statements about which types occur. Since every Delaunay
triangulation gives a point of $P_i$, the conclusions for the
infeasible type and for the types with a coordinate forced to $1$ are
proved by the linear system alone, without the nonlinear
circle-pattern theory of~\cite{Springborn2003,BobenkoSpringborn2004}.

\subsection{Loops}

The three full-dimensional types with loops ($13,17,19$) show that a
bridge of $G$ does not by itself obstruct a full-dimensional solution
space: once $\rank(A_i)=6$ (Lemma~\ref{lem:rankwithloops}), the loop's
column behaves like any other. In the three lower-dimensional types
with loops ($15,16,18$) every loop coordinate is forced to $0$, and
each of these types also has a coordinate forced to $1$, so by
Proposition~\ref{prop:degenerate}(b) it does not occur.

\subsection{Further directions}

Two questions remain. First, the gluing of
Section~\ref{sec:facets} concerns the top cells only; extending it to
the faces of all codimensions, and to the realizable types of the
contraction hierarchy that appear as the polytopes of types
$3,5,10,11$ (Proposition~\ref{prop:degenerate}(c)), would give the
full combinatorial cell structure. Its Euler characteristic can then
be compared with $\chi^{\mathrm{orb}}(\mathcal M_2)=-1/240$. Second,
whether the wedge-angle map identifies $P_i^\circ$ with the
corresponding cell of Teichm\"uller space is the nonlinear
realizability question of~\cite[Section~7]{Paper1}; the polytopes
computed here are its linear necessary conditions.

\section{Conclusion}

The rank-six phenomenon of~\cite{Paper1}, proved there for the two
generic Delaunay types, holds for all twenty trivalent types,
including the six whose triangulation has a loop. This uniformity
hides a sharp dichotomy at the level of the polytope of wedge angles:
exactly the ten realizable types have a six-dimensional polytope, nine
have lower-dimensional ones, and one type has none. The lower
dimensions have a geometric meaning: four of these polytopes belong
to realizable types with fewer edges, and the other six types occur
as no Delaunay triangulation at all. The vertices of all twenty
polytopes are half-integral. For the ten full-dimensional types,
facets occur only at Whitehead-admissible edges, a Whitehead move
identifies the facets of the two adjacent polytopes, and these
identifications are compatible with the automorphism groups, giving a
finite complex of ten orbifold cells glued along $26$ walls. Its free
unfolding grows at rate about $7.093$; the true cell decomposition of
Teichm\"uller space is a quotient of this unfolding in which the cells
around codimension-two faces close up.

\section*{Code availability}

The verification script and the Lean~4 formalization of
Section~\ref{sec:verification} are in the repository
\url{https://github.com/armandorodado/weierstrass-delaunay}, folder
\texttt{paper2-\allowbreak circle-\allowbreak pattern-\allowbreak systems},
archived at Zenodo~\cite{Code}.

\section*{Tool use}

\emph{Software.} The computations of Section~\ref{sec:verification}
were carried out in Python~3 with \texttt{numpy}, \texttt{scipy}
(the HiGHS linear-programming solver) and \texttt{networkx}, and the
formal proofs were checked with Lean~4 (version 4.15.0, core library
only). The paper was typeset in \LaTeX{} and the figure drawn with
Ti\textit{k}Z.

\emph{Artificial intelligence.} The AI assistant Claude (Anthropic)
was used in preparing this version of the paper, as follows: to
review the draft for mathematical correctness and consistency with
the companion paper~\cite{Paper1}, which led to the corrections and
the added proofs of Lemma~\ref{lem:degreetwo},
Theorem~\ref{thm:gluing} and Proposition~\ref{prop:degenerate}; to
write the independent verification script and the Lean~4
formalization of Section~\ref{sec:verification}; and to edit the text
for language, structure and \LaTeX{} typesetting. All mathematical
content, code and output were checked by the author, who takes full
responsibility for the paper.

\end{document}